\documentclass[12pt]{article}
\usepackage{amssymb}
\usepackage{mathrsfs}
\usepackage{amsthm}
\usepackage{amsmath}
\usepackage{graphicx}
\usepackage{tikz}
\usepackage{enumerate}

\usepackage{latexsym,amsmath,amssymb,amsfonts,cite}
\usepackage{eepic,color,colordvi,amscd}

\newtheorem{thm}{Theorem}

\newtheorem{clm}{Claim}

\newcommand*{\eqdef}{\stackrel{\text{\tiny{def}}}{=}}            % definition by equality
\newcommand*{\abs}[1]{\lvert #1\rvert}                           % Absolute values, cardinality
\DeclareMathOperator{\LCS}{LCS}
\DeclareMathOperator{\rev}{rev}
\DeclareMathOperator{\dist}{dist}
\newcommand*{\0}{\textsf{0}}
\newcommand*{\1}{\textsf{1}}
\newcommand*{\2}{\textsf{2}}
\newcommand*{\3}{\textsf{3}}
\newcommand*{\4}{\textsf{4}}
\newcommand*{\z}{\textsf{0}}
\newcommand*{\symt}{\textsf{t}}
\newcommand*{\dLP}{\dot{LP}}
\newcommand*{\dRP}{\dot{RP}}
\newcommand*{\dw}{\dot{w}}

\begin{document}

\title{Longest common subsequences in sets of words}
\author{Boris Bukh\footnote{Department of Mathematical Sciences, Carnegie Mellon University, Pittsburgh, PA 15213, USA. Email: bbukh@math.cmu.edu. Supported in part by U.S. taxpayers through NSF grant DMS-1201380.}\and
Jie Ma\footnote{School of Mathematical Sciences, University of Science and Technology of China, Hefei, Anhui 230026, China.}\\
}

\date{}

\maketitle

\begin{abstract}
Given a set of $t\ge k+2$ words of length $n$ over a $k$-letter alphabet, it is proved that there exists a common subsequence among two of them of length at least $\frac{n}{k}+cn^{1-1/(t-k-2)}$, for some $c>0$ depending on $k$ and $t$. This is sharp up to the value of $c$.
\end{abstract}

\section{Introduction}
A {\it word} is a sequence of symbols from some fixed finite alphabet.
For the problems in this paper only the size of the alphabet is important.
So we will use $[k]\eqdef \{1,2,\ldots, k\}$ for a canonical $k$-letter alphabet.
The family of all words of length $n$ over a $k$-letter alphabet is thus denoted by $[k]^n$.
For a word $w$, a {\it subsequence} is any word obtained by deleting zero or more symbols from $w$.
By a {\it subword} of $w$, we mean a subsequence of $w$ consisting of consecutive symbols.
For example, $\1\3\3\4$ is a subsequence but not a subword of $\1\2\3\4\1\2\3\4$.
A {\it common subsequence} of $w$ and $w'$ is a word that is a subsequence of both $w$ and $w'$.

A general principle asserts that every sufficiently large collection of objects necessarily contains a pair of similar objects.
In this paper, we treat the case when the objects are words, and similarity is measured by length of a common subsequence.
We use $\LCS(w,w')$ to denote the length of the longest common subsequence of words $w$ and $w'$.
For a set $\mathcal{W}$ of words, let $\LCS(\mathcal{W})\eqdef \max \LCS(w,w')$ where the maximum is taken over all pairs $\{w,w'\}$ in $\mathcal{W}$.
We also allow $\mathcal{W}$ to be a multiset, so $\mathcal{W}$ might contain some elements multiple times.

For an integer $t\ge 2$ and a family $\mathcal{F}$ of words,
let
  \[
    \LCS(t,\mathcal{F})\eqdef \min_{\mathcal{W}\in \mathcal{F}^t} \LCS(\mathcal{W}).
  \]

A \emph{permutation} of length $k$ is a word over $[k]$ in which every symbol appears exactly once. Let $\mathcal{P}_k$ be the set of all permutations of length $k$.
Much of the inspiration
for our work comes from the results of Beame--Huynh-Ngoc \cite{BH}, Beame--Blais--Huynh-Ngoc\cite{BBH} and Bukh--Zhou \cite{BZ} that can
be summarized as
\begin{align*}
  \LCS(3, \mathcal{P}_k)&=k^{1/3}+O(1),\\
  \LCS(4, \mathcal{P}_k)&=k^{1/3}+O(1),\\
  1.001k^{1/3}\leq \LCS(t, \mathcal{P}_k)&\leq 4k^{1/3}+O(k^{7/40})&&\text{for }5\leq t\leq k^{1/3}.
\end{align*}
The problem of bounding $\LCS(t, \mathcal{P}_k)$ is closely related to the longest twin problem
of Axenovich, Person and Puzynina \cite{APP}. Here, two subsequences $w_1,w_2$ of the same word $w$ are {\it twins} if they are equal as words
but the sets of positions of symbols from $w$ retained in $w_1$ and in $w_2$ are disjoint.
It was shown in \cite{BZ} that if $\LCS(t, \mathcal{P}_k)$ is small for $t\geq 2k$, then
there are words that contain no long twins, and that a converse (which is more technical to state) also holds.

In this paper, we consider $\LCS(t,[k]^n)$.
We let $w^m$ be the concatenation of $m$ copies of $w$, e.g., $(\3\4\3)^2=\3\4\3\3\4\3$.
For $t\leq k$ we have $\LCS(t,[k]^n)=0$ as the family
$\{\1^n,\2^n,\dotsc,\symt^n\}$ shows. For $t=k+1$ we have\footnote{For readability, we omit the floor and ceiling signs throughout the paper.}
\[
  \LCS(k+1, [k]^n)=\tfrac{n}{k}
.\]
The upper bound is attained by the family $\{\1^n, \2^n, \ldots, \textsf{k}^n, (\1\2...\textsf{k})^{\frac{n}{k}}\}$ of $k+1$ words.
The lower bound is a consequence of two simple facts: the most popular letter in a word occurs at least $\frac{n}{k}$ times,
and the most popular letter is the same in two out of $k+1$ words.
The following is our main theorem, which determines the asymptotic magnitude of $\LCS(t, [k]^n)$ for all $t\ge k+2$.
\begin{thm}
\label{thm:main}
For nonnegative integers $k, r$ and $n$ such that $k\ge 2$ and $n\ge k(10r)^{9r}$, there exists $c=\Theta(r^{-9}k^{1/r-2})$ such that $$\LCS(r+k+2, [k]^n)\ge \frac{n}{k}+cn^{1-\frac1r}.$$
\end{thm}

This theorem is sharp up to the value of $c$. For $0\leq i\leq r$, put
\begin{align*}
  m_i&\eqdef (n/k)^{i/r},\\
\intertext{and define words}
  w_i&\eqdef (\1^{m_i}\2^{m_i}\ldots\textsf{k}^{m_i})^{\frac{n}{km_i}},\\
  \rev w_i&\eqdef (\textsf{k}^{m_i}\ldots\2^{m_i}\1^{m_i})^{\frac{n}{km_i}}.
\end{align*}
Note that $\rev w_i$ is the word obtained by reversing the symbols of $w_i$. We claim that for the family
\begin{equation}
\label{def:W}
  \mathcal{W}\eqdef \left\{w_0,w_1,\ldots, w_{r},\rev w_r, \1^n,\2^n, \ldots, \textsf{k}^n\right\}
\end{equation}
we have
\[
  \LCS(\mathcal{W})\le \frac{n}{k}+k^{1/r} n^{1-1/r}.
\]
Indeed, $\LCS(w_i,\textsf{j}^n)=\frac{n}{k}$ is clear, and every common subsequence of $w_r$ and $\rev w_r$ is of the form $\textsf{i}^m$, and so is of length at most $n/k$.
Hence, it suffices to bound $\LCS(w_i,w_j)$ and $\LCS(w_i,\rev w_j)$ for $i<j$ (though we need a bound
on $\LCS(w_i,\rev w_j)$ only for $j=r$). The two cases are similar:
Any common subsequence of $w_i$ and $w_j$ (or $\rev w_j$) must be of the form $k_1^{p_1}k_2^{p_2}\ldots k_s^{p_s}$ for some $s\le \frac{n}{m_j}$,
where $k_l\in [k]$ for $l=1,2,\ldots,s$. Since each subsequence $k_l^{p_l}$ of $w_i$ spans a subword of length at least $(\lceil\frac{p_l}{m_i}\rceil-1)km_i\ge (p_l-m_i)k$ in $w_i$, it follows that $n\ge \sum_{l=1}^s(p_l-m_i)k$,
implying that $\LCS(w_i,w_j)=\max \sum_l p_l\le \frac{n}{k}+\frac{m_i}{m_j}n\le \frac{n}{k}+k^{1/r} n^{1-1/r}$.

\medskip

A word $w\in [k]^n$ is called {\it balanced}
if it contains the same number of $\textsf{i}$'s as $\textsf{j}$'s for any $\textsf{i}, \textsf{j}\in [k]$.
Let $B_k^n$ be the family containing all balanced words in $[k]^n$.
As we shall see, the assertion of Theorem~\ref{thm:main} reduces to the following result on family~$B_k^n$.
\begin{thm}
\label{thm:balanced}
For nonnegative integers $k, r$ and $n$ such that $k\ge 2$ and $n\ge k(10r)^{9r}/2$, there exists $c'=\Theta(r^{-9}k^{1/r-1})$ such that $$\LCS(r+2, B_k^n)\ge \frac{n}{k}+c'n^{1-\frac1r}.$$
\end{thm}
For large $n$, Theorem~\ref{thm:balanced} is sharp up to the value of $c'$, as witnessed by the family of words $\{w_0,w_1,\ldots, w_r,\rev w_r\}$,
where $w_i$'s are as in \eqref{def:W}. For values of $n$ that are comparable to $k$, the rates of growth as $k\to \infty$
for $\LCS(3,B_k^n)$ and for $\LCS(4,B_k^n)$ have been determined in \cite[Theorems 11 and 12]{BZ}.

\medskip

The rest of the paper is organized as follows.
In the next section, we reduce Theorems \ref{thm:main} and~\ref{thm:balanced} to the $\LCS$ for balanced binary words (see Theorem~\ref{thm:binary}).
The proof of Theorem~\ref{thm:binary} will be completed in section~\ref{sec:proof}. The last section contains a couple of open problems.
In this paper, we do not attempt to optimize the constants, and instead aim for simpler presentation.

\section{Reductions}
In this section we deduce Theorems \ref{thm:main} and \ref{thm:balanced} from the following special case $k=2$ of Theorem~\ref{thm:balanced}, which we state separately.
\begin{thm}
\label{thm:binary}
Let $r$ and $n$ be nonnegative integers such that $n\ge (10r)^{9r}$. For any set $\mathcal{W}$ of $r+2$ balanced words in $\{\0,\1\}^n$,
we have $$\LCS(\mathcal{W})\ge \frac{n}{2}+ \Omega\left(r^{-9}\right)\cdot n^{1-1/r}.$$
\end{thm}

\noindent{\it Proof of Theorem~\ref{thm:balanced}.} ({\it Assume that Theorem~\ref{thm:binary} holds.})
Consider a multiset $\mathcal{W}$ of arbitrary $r+2$ balanced words from $[k]^n$ and let $$\mathcal{W}'\eqdef \{\text{the subsequence of } w \text{ consisting of all } \1 \text{ 's and } \2 \text{ 's for every } w\in \mathcal{W}\}.$$
Then $\mathcal{W}'$ is a multiset of $r+2$ balanced words from $\{\1,\2\}^{n'}$ for $n'=\frac{2n}{k}\geq (10r)^{9r}$.
By Theorem \ref{thm:binary}, we have $\LCS(\mathcal{W})\ge \LCS(\mathcal{W}')\ge \frac{n'}{2}+\Omega(r^{-9}) (n')^{1-1/r}=\frac{n}{k}+\Omega(r^{-9}k^{1/r-1})\cdot n^{1-1/r}.$\qed

\bigskip

\noindent{\it Proof of Theorem~\ref{thm:main}.} ({\it Assume that Theorem~\ref{thm:balanced} holds.})
Let $c$ be a small constant. Consider an arbitrary set of $r+k+2$ words from $[k]^n$.
Call $w\in [k]^n$ {\it unhinged} if some letter occurs in $w$ at least $\frac{n}{k}+cn^{1-1/r}$ times and {\it hinged} otherwise.
Observe that if there are $k+1$ unhinged words in the set, then some two of them have $\LCS$ of length at least $\frac{n}{k}+cn^{1-1/r}$.
Thus we may assume that there are at least $r+2$ hinged words.
Since each hinged word of length $n$ contains a subsequence that is a balanced word of length $n-k^2cn^{1-1/r}$,
by Theorem~\ref{thm:balanced}, some two hinged words have $\LCS$ of length at least
$$\frac{n-k^2cn^{1-\frac1r}}{k}+c'\left(n-k^2cn^{1-\frac1r}\right)^{1-\frac1r}\ge \frac{n}{k}+cn^{1-\frac1r},$$ provided $c=\Theta\left(\frac{c'}{k}\right)$.
This proves Theorem~\ref{thm:main}. \qed

\section{The proof of Theorem~\ref{thm:binary}}\label{sec:proof}
Throughout this proof, let $\alpha\eqdef 10^{-6}r^{-9}$ and $\beta\eqdef\frac{1}{40000}r^{-6}$, and let $\mathcal{W}$ be a set consisting of arbitrary balanced words $w^{(1)},w^{(2)},\ldots,w^{(r+2)}$ in $\{\0,\1\}^n$.
Moreover, we assume that $r\ge 2$, as it is easy to see that
$\LCS(\mathcal{W})\ge \frac{n}{2}$ when $r=0$ and $\LCS(\mathcal{W})\ge \frac{n}{2}+1$ when $r=1$.

We give a very brief outline of the proof before proceeding.
A crucial idea is to consider the scale on which $\0$'s and $\1$'s alternate in a word.
For example, in word $(\0\1)^{n/2}$ alternation (between $\0$'s and $\1$'s) happens on scale $\Theta(1)$,
whereas in word $\0^{n/2}\1^{n/2}$ the alternation scale is about $\Theta(n)$.
The proof will first find two words, say $w^{(1)}$ and $w^{(2)}$, of ``comparable" alternation scale,
and then show, in effect, that $\LCS(w^{(1)},w^{(2)})$ is large.

We shall think of words as made of a sequence of
\emph{distinguishable} $\0$'s and $\1$'s. That means that if we say ``let $z$ be a $\0$
in word $w$'', then the variable $z$ refers to a particular $\0$. For example,
if $z$ is the $3$'rd $\0$ in the word $\1\1\0\0\1\0\1\1\0\1$, and $w'$ is the word
obtained from $w$ by removing the $1$'st and $4$'th zeros, namely $w'=\1\1\0\1\0\1\1\1$,
then $z$ becomes the $2$'nd zero in $w'$.

For two symbols $a,b$ from $w$ such that $a$ is to the left of $b$,
we denote by $w[a,b]$ the subword of $w$ starting from $a$ and ending with $b$.
We also use $w(a,b)$ to denote the subword of $w$ obtained from $w[a,b]$ by deleting $a$ and $b$.
The notations $w[a,b)$ and $w(a,b]$ are defined similarly.

If $z$ is a $\0$ in a word $w$, its \emph{position}, denoted $P_w(z)$, is the number of $\1$'s to the
left of $z$. When the word $w$ is clear from the context, we will drop the subscript of $P_w(z)$ and write simply $P(z)$.
Note that several $\0$'s might have the same position. If $z$ is the $j$'th
$\0$ in $w$, we say that its \emph{expected position} is $j$, because in a random word the expected value of $P(z)$ is $j$.
We say that a $\0$ is \emph{good}
in $w$ if its position differs from its expected position by at most $\alpha n^{1-1/r}$.
If a $\0$ is not good, then its position is either to the left or to the right of
its expected position. In these cases we call such a $\0$ \emph{left-bad} and \emph{right-bad}
respectively. The following claim will be used frequently.

\smallskip

\begin{clm}\label{clm1}
If a subword of $w^{(i)}$ contains $N$ $\1$'s, then it contains at most $N+2\alpha n^{1-1/r}$ good $\0$'s.
\end{clm}
\begin{proof}
Let $z_l$ and $z_r$ be the leftmost and the rightmost good $\0$'s in the subword. By definition, we have $P(z_r)-P(z_l)\leq N$.
From the goodness of $z_r$, we see that its expected position differs from $P(z_r)$ by at most $\alpha n^{1-1/r}$; similarly it holds for $z_l$.
Therefore, the expected positions of $z_r$ and $z_l$ differ by at most $N+2\alpha n^{1-1/r}$, implying this claim.
\end{proof}

\smallskip

We introduce a concept closely related to the alternation scale described in the outline.
A subword is called a \emph{$\0$-rich interval of length $L$} if it contains exactly
$L$ good $\0$'s and no more than $L/10$ $\1$'s.
The \emph{type} of a good $\0$ is the
largest integer $t$ such that this $\0$ is contained in a $\0$-rich interval of length exactly $n^{t/r}$.
%We will omit rounding signs for clarity of presentation,
%thus throughout the proof we regard $n^{t/r}$ and other parameters as integers (if necessary).
Note that a type of a good $\0$ is well-defined since every good zero is contained in a
$\0$-rich interval of length~$1$. Also note that a type cannot be $r-1$.
Indeed, if there existed a $\0$-rich interval of length $n^{1-1/r}$, then Claim~\ref{clm1} would imply
that $n^{1-1/r}\le n^{1-1/r}/10+2\alpha n^{1-1/r}$, a contradiction.
We define a \emph{type of a bad $\0$} to be either \textsf{left-bad} or \textsf{right-bad}.
Thus a type of each $\0$ is an element of $\{0,1,\dotsc,r-2,\textsf{left-bad},\textsf{right-bad}\}$.

To be able to refer to individual $\0$'s, we define $\z_j^{(i)}$ as the $j$'th $\0$ in word $w^{(i)}$. As
our proof does not treat $\0$'s and $\1$'s symmetrically, we do not need a similar notation to refer
to individual $\1$'s.

Fix an integer $j$ and consider $\z_j^{(1)},\z_j^{(2)}, \dotsc,\z_j^{(r+2)}$.
We may assume that at most one of these zeros is left-bad, and at most one of them is right-bad.
Suppose, on the contrary, that both $\z_j^{(1)}$ and $\z_j^{(2)}$ are left-bad. Then
we can obtain a common subsequence of $w^{(1)}$ and $w^{(2)}$ with length at least $n/2+\alpha n^{1-1/r}$ by
matching up the first $j$ $\0$'s and then $\1$'s to the right of $\z_j^{(1)}$ and $\z_j^{(2)}$.
Hence, in this case $\LCS(\mathcal{W})\ge \LCS(w^{(1)},w^{(2)})\geq n/2+\alpha n^{1-1/r}=n/2+\Omega(r^{-9})\cdot n^{1-1/r}$. The case of two right-bad $\0$'s is similar.

Hence, for any integer $j$, two of $\z_j^{(1)},\z_j^{(2)},\dotsc, \z_j^{(r+2)}$ are of the same type,
and that type is one of $0,1,\dotsc, r-2$.
By the pigeonhole principle, there are two words, say $w^{(1)}$ and $w^{(2)}$, and some $t\in \{0,1,\ldots, r-2\}$ such that the set
\[
  \mathcal{T}\eqdef \{ j  : \text{ both } \z_j^{(1)} \text{ and } \z_j^{(2)} \text{ have type } t\}
\]
has size at least $\frac{n/2}{\binom{r+2}{2}(r-1)}\ge \frac{n}{2r^3}$. We will show that
$w^{(1)}$ and $w^{(2)}$ contain a common subsequence of length $n/2+\Omega(n^{1-1/r})$.

\medskip

We partition each of $w^{(1)}$ and $w^{(2)}$ into {\it blocks} that contain exactly $\beta n^{1-1/r}$ many $\1$'s.
To be more precise, for each $i\in \{1,2\}$, the $k$'th block (denoted by $B_k^{(i)}$) of word $w^{(i)}$ is defined to be the subword $w^{(i)}[a_{k-1},a_k)$, where $a_k$ denotes the $(k\cdot \beta n^{1-1/r}+1)$'th $\1$ in word $w^{(i)}$.

For each $i\in \{1,2\}$ and each $j\in \mathcal{T}$, choose a $\0$-rich interval of length $n^{t/r}$ containing $\z_j^{(i)}$ and call this interval $I_j^{(i)}$. By shrinking $I_j^{(i)}$ if necessary, we may assume that both leftmost and rightmost symbols in $I_j^{(i)}$ are good $\0$'s.
An integer $j\in \mathcal{T}$ is {\it consistent} if $I_j^{(1)}\subset B_{k}^{(1)}$ and $I_j^{(2)}\subset B_{k}^{(2)}$ for some $k$.
Let $\mathcal{S}=\{j\in \mathcal{T} : j\text { is consistent}\}$.

\smallskip

\begin{clm}\label{clm2}
 $|\mathcal{S}|\ge \frac{n}{4r^3}$.
\end{clm}
\begin{proof}
For each $i\in \{1,2\}$, let $L_{k}^{(i)}$ be the subword of $w^{(i)}$ spanning the last $2\alpha n^{1-1/r}$ $\1$'s in the block $B_{k}^{(i)}$ and the first $2\alpha n^{1-1/r}$ $\1$'s in the block $B_{k+1}^{(i)}$.
By Claim~\ref{clm1}, we see that $L_{k}^{(i)}$ contains at most $6\alpha n^{1-1/r}$ good $\0$'s,
and hence the set
$L^{(i)}\eqdef \{ \text{all good }\0\text{'s contained in }\cup_{k} L_{k}^{(i)}\}$
is of size at most $6\alpha n^{1-1/r}\cdot \frac{n/2}{\beta n^{1-1/r}}=\frac{3\alpha n}{\beta}$.
Let $$\mathcal{T}'\eqdef \{j\in \mathcal{T}: \0_j^{(i)}\notin L^{(i)} \text{ for each } i=1,2\}.$$
It is clear that $|\mathcal{T}'|\ge |\mathcal{T}|-|L^{(1)}|-|L^{(2)}|\ge \frac{n}{2r^3}-\frac{6\alpha n}{\beta}\ge \frac{n}{4r^3}$.

Now it suffices to show that $\mathcal{T}'\subseteq \mathcal{S}$.
Consider an arbitrary integer $j\in \mathcal{T}'$.
Assume that $\z_j^{(1)}\in B_k^{(1)}$ for some $k$. By the definition of $\mathcal{T}'$, it holds that
$$(k-1)\beta n^{1-1/r}+2\alpha n^{1-1/r}< P(\z_j^{(1)})\le k\beta n^{1-1/r}-2\alpha n^{1-1/r}.$$
As $t\le r-2$ and $n\ge (10r)^{9r}$, the $\0$-rich interval $I_j^{(1)}$ has at most $n^{t/r}/10\le 2\alpha n^{1-1/r}$ $\1$'s, implying that $I_j^{(1)}\subset B_k^{(1)}$.
By the goodness of $\z_j^{(1)}$ and $\z_j^{(2)}$, we obtain that $|P(\z_j^{(1)})-P(\z_j^{(2)})|\le 2\alpha n^{1-1/r}$,
which implies that $\z_j^{(2)}\in B_k^{(2)}$.
By the definition of $\mathcal{T}'$ again, in fact we have
$(k-1)\beta n^{1-1/r}+2\alpha n^{1-1/r}< P(\z_j^{(2)})\le k\beta n^{1-1/r}-2\alpha n^{1-1/r}.$
Repeating the same argument, we see $I_j^{(2)}\subset B_k^{(2)}$.
So $j$ is consistent and hence $j\in \mathcal{S}$, finishing the proof of Claim~\ref{clm2}.
\end{proof}
\smallskip

With slight abuse of notation, let $\mathcal{S}\cap B_k^{(i)}\eqdef \{\z_j^{(i)}\in B_k^{(i)}: j\in \mathcal{S}\}$.
Clearly, $\mathcal{S}\cap B_k^{(1)}$ and $\mathcal{S}\cap B_k^{(2)}$ are of the same size, say $s_k$. Then $s_k$ satisfy
\begin{align}
\label{equ:sk}
0\le s_k\le (\beta+2\alpha)n^{1-1/r} \text{~~ and~~} \sum_{k} s_k= |\mathcal{S}|\ge \frac{n}{4r^3},
\end{align}
where the first inequality follows by Claim~\ref{clm1}.
For fixed $k$ and $i\in \{1,2\}$, consider the family of all $\0$-rich intervals $I_j^{(i)}$
that belong to $B_k^{(i)}$ as $j$ ranges over $\mathcal{S}$.
It is clear that the union of $I_j^{(i)}$'s from this family contains all $\0$'s in $\mathcal{S}\cap B_k^{(i)}$.
By the Vitali covering lemma, there is a subfamily, denoted by $\mathcal{I}_k^{(i)}$, consisting of pairwise disjoint intervals $I_j^{(i)}$
whose union contains at least one third of the $\0$'s in $\mathcal{S}\cap B_k^{(i)}$.
Since each $I_j^{(i)}$ contains at most $n^{t/r}$ $\0$'s from $\mathcal{S}\cap B_k^{(i)}$, we derive
\begin{align}
\label{equ:Iik}
|\mathcal{I}_k^{(i)}|\ge \frac{s_k}{3n^{t/r}}.
\end{align}
Let $\mathcal{I}^{(i)}\eqdef\cup_{k} \mathcal{I}_k^{(i)}$. The intervals in $\mathcal{I}^{(i)}$ are disjoint, for
intervals in $\mathcal{I}_k^{(i)}$ and $\mathcal{I}_{k'}^{(i)}$ for $k\neq k'$ are contained in non-overlapping blocks.

\smallskip

We shall pick an integer $Q$ in the interval $(-\beta n^{1-1/r}, \beta n^{1-1/r})$ uniformly at random,
and define words $\dot{w}^{(1)}$ and $\dot{w}^{(2)}$ as follows.
If $Q\ge 0$, let $\dot{w}^{(1)}\eqdef w^{(1)}$ and $\dot{w}^{(2)}$ be obtained from $w^{(2)}$ by removing the first $Q$ $\1$'s;
otherwise, let $\dot{w}^{(2)}\eqdef w^{(2)}$ and $\dot{w}^{(1)}$ be obtained from $w^{(1)}$ by removing the first $-Q$ $\1$'s.
%Note that any interval $I\in \mathcal{I}_i$ remains an interval in $w_i'$, as the changes are only made in two $1$'st blocks.
%And we will use $\dot{LP}(I)$ and $\dot{RP}(I)$ to denote the left-position and the right-position of an interval $I$ in $\dot{w}^{(i)}$.

For an interval $I\in \mathcal{I}^{(i)}$, its {\it left-position} (resp.\ {\it right-position}) in $w^{(i)}$
is the position of the leftmost (resp.\ rightmost) good $\0$ in $w$. We denote left- and right-positions by $LP(I)$ and $RP(I)$.
We define the left- and right-positions of an interval in $\dw^{(i)}$ similarly, and denote them by $\dLP(I)$ and $\dRP(I)$.
We note that for intervals $I_1\in \mathcal{I}^{(1)}$ and $I_2\in\mathcal{I}^{(2)}$
\begin{equation}\label{equ:positions}
  LP(I_2)-LP(I_1)-Q=\dLP(I_2)-\dLP(I_1).
\end{equation}

We say that two intervals $I_1\in \mathcal{I}^{(1)}$ and $I_2\in \mathcal{I}^{(2)}$ are {\it close} or $(I_1,I_2)$ is a {\it close pair}, if
\begin{align*}
|\dLP(I_2)-\dLP(I_1)|\le \frac{1}{20}n^{t/r}.
\end{align*}
Suppose that intervals $I_1\in \mathcal{I}^{(1)}$ and $I_2\in \mathcal{I}^{(2)}$ are close, then as $0\le \dRP(I_i)-\dLP(I_i)\le \frac{1}{10}n^{t/r}$, we also have
\begin{align}
\label{equ:RP}
|\dRP(I_2)-\dRP(I_1)|\le \frac{3}{20}n^{t/r}.
\end{align}

\smallskip

\begin{clm}\label{clm3}
Each interval in $\mathcal{I}^{(1)}$ is close to at most $n^{1/r}$ intervals in $\mathcal{I}^{(2)}$.
Similarly, each interval in $\mathcal{I}^{(2)}$ is close to at most $n^{1/r}$ intervals in $\mathcal{I}^{(1)}$.
\end{clm}

\begin{proof}
Suppose, on the contrary, that an interval $I\in \mathcal{I}^{(1)}$ is close to $J_1,J_2,\dotsc,J_{d}\in \mathcal{I}^{(2)}$
with $\dLP(J_1)<\dLP(J_2)<\ldots<\dLP(J_{d})$, where $d\eqdef n^{1/r}+1$.
Let $J$ be the subword of $w^{(2)}$ starting from the leftmost good $\0$ of $J_1$ and ending with the leftmost good $\0$ of $J_d$.
By the closeness of $(I,J_1)$ and of $(I,J_d)$, we have $|\dLP(J_1)-\dLP(J_d)|\le \frac{1}{10}n^{t/r}$,
which implies that $J$ has at most $\frac{1}{10}n^{t/r}\le \frac{1}{10}n^{(t+1)/r}$ $\1$'s.
Since $J$ also contains at least $(d-1)\cdot n^{t/r}=n^{(t+1)/r}$ good $\0$'s,
every $\0$ in $J$ is contained in a $\0$-rich interval of length $n^{(t+1)/r}$.
Hence the type of any $\0$ in $J_1$ is at least $t+1$.
Yet from the construction of $\mathcal{I}^{(2)}$, it is evident that $J_1$ contains at least one $\0$ of type $t$.
This contradiction finishes the proof of Claim~\ref{clm3}.
\end{proof}

Some intervals in the first block, i.e., those in $\mathcal{I}_1^{(1)}\cup \mathcal{I}_1^{(2)}$, might be destroyed
in the passage from $w^{(1)}$ and $w^{(2)}$ to their dotted counterparts. So let $k\geq 2$ and
consider two arbitrary intervals $I_1\in \mathcal{I}_k^{(1)}$ and $I_2\in \mathcal{I}_k^{(2)}$.
In view of \eqref{equ:positions}, $I_1$ and $I_2$ are close if and only if
\begin{align}
\label{equ:close}
|LP(I_2)-LP(I_1)-Q|\le \frac{1}{20}n^{t/r}.
\end{align}
Since $I_1$ and $I_2$ are in the same block, there exists
an integer $q\in (-\beta n^{1-1/r},\beta n^{1-1/r})$ such that $LP(I_2)=LP(I_1)+q$.
Therefore there are at least $\frac{1}{20}n^{t/r}$ choices of $Q$'s for which \eqref{equ:close} holds,
namely $Q$ can be any integer in $[q-\frac{1}{20}n^{t/r},q+\frac{1}{20}n^{t/r}]\cap (-\beta n^{1-1/r},\beta n^{1-1/r})$.
This shows that the probability that $I_1$ and $I_2$ are close is at least
\begin{align}
\label{equ:prob}
p\eqdef\frac{\frac{1}{20}n^{t/r}}{2\beta n^{1-1/r}}=\frac{1}{40\beta}n^{1/r+t/r-1}.
\end{align}

Let $E\subset \mathcal{I}^{(1)}\times \mathcal{I}^{(2)}$ be the set of close pairs $(I_1,I_2)$.
Then the expectation of $|E|$ is at least
$p\cdot \left(\sum_{k\ge 2}|\mathcal{I}_k^{(1)}| |\mathcal{I}_k^{(2)}|\right)$.
There must exist some $Q\in (-\beta n^{1-1/r}, \beta n^{1-1/r})$ such that the size of $E$ is at least its expectation.
Fix such a $Q$. Note that this also fixes $\dot{w}^{(1)}, \dot{w}^{(2)}$ and the set $E$.
By \eqref{equ:sk}, \eqref{equ:Iik}, \eqref{equ:prob} and the Cauchy--Schwarz inequality, we derive
\begin{align}
\label{equ:E}
|E|\ge p\cdot \left(\sum_{k\ge 2}|\mathcal{I}_k^{(1)}| |\mathcal{I}_k^{(2)}|\right)\ge \frac{n^{1-t/r}}{5000r^6},
\end{align}
since the summation is over at most $\frac{n/2}{\beta n^{1-1/r}}=\frac{n^{1/r}}{2\beta}$ terms.

\smallskip

\begin{clm}\label{clm4}
There exist $\frac{|E|}{2n^{1/r}}$ close pairs $(I_i,J_i)$ in $E$
such that
\begin{equation}\label{equ:matching}
\begin{aligned}
\dLP(I_1)&<\dLP(I_2)<\ldots <\dLP(I_{|E|/2n^{1/r}}),\\
\dLP(J_1)&<\dLP(J_2)<\ldots <\dLP(J_{|E|/2n^{1/r}}).
\end{aligned}
\end{equation}
\end{clm}

\begin{proof} We can view $E$ as the edge set of a bipartite graph $G$ with bipartition $(\mathcal{I}^{(1)},\mathcal{I}^{(2)})$.
We desire to find a large matching $I_1J_1,I_2J_2,\dotsc$ satisfying \eqref{equ:matching}.
Identify $I\in \mathcal{I}^{(1)}$ with the point $(\dLP(I),0)$ in the Euclidean plane, and identify $J\in\mathcal{I}^{(2)}$ with
the point $(\dLP(J),1)$. Edges will be represented by line segments.
Among all the matchings of maximum size, pick one that minimizes
the total Euclidean length of edges.

We claim that this matching satisfies \eqref{equ:matching}. Suppose, on the contrary, that $\dLP(I_l)<\dLP(I_m)$ and
$\dLP(J_l)>\dLP(J_m)$. Then the line segments $\overline{I_lJ_l}$ and $\overline{I_mJ_m}$ cross. That implies
that line segments $\overline{I_lJ_m}$ and $\overline{I_mJ_l}$ are both shorter than $\max\bigl(\dist(I_l,J_l),\dist(I_l,J_l)\bigr)$,
and so $I_lJ_m,I_mJ_l\in E$. This contradicts the choice of the matching, since replacing
edges $I_lJ_l$ and $I_mJ_m$ with $I_lJ_m$ and $I_mJ_l$ decreases the total length of the matching as the following picture demonstrates.

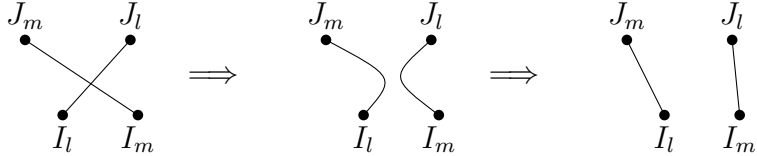
\begin{center}
\begin{tikzpicture}
\begin{scope}
\coordinate (Il) at (0,0);
\coordinate (Im) at (1,0);
\coordinate (Jl) at (0.9,1);
\coordinate (Jm) at (-0.5,1);
\fill (Il) circle (2pt) node [below] {$I_l$};;
\fill (Im) circle (2pt) node [below] {$I_m$};;
\fill (Jl) circle (2pt) node [above] {$J_l$};;
\fill (Jm) circle (2pt) node [above] {$J_m$};;
\draw (Il)--(Jl);
\draw (Im)--(Jm);
\end{scope}
\begin{scope}[xshift=4cm]
\coordinate (Il) at (0,0);
\coordinate (Im) at (1,0);
\coordinate (Jl) at (0.9,1);
\coordinate (Jm) at (-0.5,1);
\fill (Il) circle (2pt) node [below] {$I_l$};;
\fill (Im) circle (2pt) node [below] {$I_m$};;
\fill (Jl) circle (2pt) node [above] {$J_l$};;
\fill (Jm) circle (2pt) node [above] {$J_m$};;
%\draw (Il) -- (0.4,0.5) -- (Jm);
%\draw (Im) -- (0.5,0.5) -- (Jl);
\draw (Il) .. controls (0.44,0.5) .. (Jm);
\draw (Im) .. controls (0.34,0.5) .. (Jl);
\end{scope}
\begin{scope}[xshift=8cm]
\coordinate (Il) at (0,0);
\coordinate (Im) at (1,0);
\coordinate (Jl) at (0.9,1);
\coordinate (Jm) at (-0.5,1);
\fill (Il) circle (2pt) node [below] {$I_l$};;
\fill (Im) circle (2pt) node [below] {$I_m$};;
\fill (Jl) circle (2pt) node [above] {$J_l$};;
\fill (Jm) circle (2pt) node [above] {$J_m$};;
\draw (Il) -- (Jm);
\draw (Im) -- (Jl);
\end{scope}
\node at (2,0.5) {$\implies$};
\node at (6,0.5) {$\implies$};
\end{tikzpicture}
\end{center}

The bound on the size of the matching follows from Claim~\ref{clm3}, which in the present language asserts that
the maximum degree in $E$ is at most $n^{1/r}$. Indeed, for any matching $M$ of size less than $\abs{E}/2n^{1/r}$
there is an $e\in E$ not adjacent to any edge of $M$. Hence, a maximal matching has at least $\abs{E}/2n^{1/r}$
edges.
\end{proof}

\smallskip

Finally, using the close pairs of Claim~\ref{clm4}, we find a long common subsequence of $w^{(1)}$ and $w^{(2)}$.
For convenience write $\lambda\eqdef |E|/2n^{1/r}$ and without loss assume that $Q\ge 0$.
For $1\le i\le \lambda-1$, let $A_i$ be the subword of $\dot{w}^{(1)}$ between the intervals $I_i$ and $I_{i+1}$ and
$B_i$ be the subword of $\dot{w}^{(2)}$ between the intervals $J_i$ and $J_{i+1}$.
In addition, let $A_0$ be the subword of $\dot{w}^{(1)}$ before the interval $I_1$, and $A_{\lambda}$ be the subword of $\dot{w}^{(1)}$ after the interval $I_{\lambda}$;
the definitions of $B_0$ and $B_{\lambda}$ are similar.
Let us consider the common subsequence $w$ of $\dot{w}^{(1)}$ and $\dot{w}^{(2)}$,
which consists of the common $\0$'s of $I_i$ and $J_i$ and the common $\1$'s of $A_i$ and $B_i$ for all $0\le i\le \lambda$.
By \eqref{equ:RP} and \eqref{equ:close}, we have
$$ |\dot{RP}(J_i)-\dot{RP}(I_i)|\le \frac{3}{20}n^{t/r}  \text{~~and~~} |\dot{LP}(J_{i+1})-\dot{LP}(I_{i+1})|\le \frac{1}{20}n^{t/r},$$
which shows that for each $i<\lambda$, the counts of $\1$'s in $A_i$ and in $B_i$ differ by at most $\frac{1}{5}n^{t/r}$.
Also note that each $I_i$ contains at most $\frac{1}{10}n^{t/r}$ $\1$'s,
thus the number of $\1$'s in $\dot{w}^{(2)}$ but not in $w$ is at most $\frac{\lambda}{2}n^{t/r}$.
By \eqref{equ:E} as well as the facts that $\lambda=|E|/2n^{1/r}$ and $|Q|< \beta n^{1-1/r}$, we derive that
\begin{align*}
\LCS(\mathcal{W})&\ge \LCS(w^{(1)},w^{(2)})\ge \LCS(\dot{w}^{(1)},\dot{w}^{(2)})\ge |w|\ge \lambda\cdot n^{t/r}+\left(\frac{n}{2}-|Q|-\frac{\lambda}{2}\cdot n^{t/r}\right)\\
&\ge \frac{n}{2}-\beta n^{1-1/r}+\frac{n^{1-1/r}}{20000r^6}=\frac{n}{2}+\Omega(r^{-6})\cdot n^{1-1/r}.
\end{align*}
This completes the proof of Theorem~\ref{thm:binary}.

\section{Two problems}
In this paper we proved that $\LCS(r+k+2, [k]^n)= \frac{n}{k}+\Theta_{r,k}(n^{1-1/r})$. It is possible that
the coefficient in the big-theta notation need not depend on~$r$, but we have been unable to prove so.
In particular, what is the smallest $r$ such that $\LCS(r+k+2, [k]^n)\geq 1.01\frac{n}{k}$? Is it
asymptotic to $\Theta(\log n)$?

Another worthy problem is the length of the longest common subsequence between two random words. A superadditivity
argument shows that the expected length of such a subsequence is asymptotic to $\gamma_k n$ for some constant $\gamma_k$.
Kiwi--Loebl--Matou\v{s}ek \cite{KLM} proved that $\gamma_k\sqrt{k}\to 2$ as $k\to\infty$, but the value
of $\gamma_k$ is not known for any $k\geq 2$ (including the case $k=4$ that is natural for the problem
of DNA comparison).

These two problems are connected. Azuma's inequality implies that $\LCS(w,w')$ for random $w,w'\in [k]^n$
is concentrated in an interval of length $\sqrt{n}$ with sub-Gaussian tails.
It thus follows that for any $\epsilon>0$ one can find a family $\mathcal{F}$ of exponentially many words from $[k]^n$
such that $\LCS(\mathcal{F})\leq\bigl(\gamma_k+\epsilon+o(1)\bigr) \cdot n$.

\end{document}